\documentstyle[12pt]{article}
\begin{document}
\begin{center}
{\bf QUANTUM ALGEBRAIC TORI}\\
 A.N.PANOV
\end{center}
{\bf Introduction.} The notion of quantum torus appears in
mathematics due to the following physical consideration.
Let the operators $u,v$ satisfy the Heisenberg commutatin relation
 $[u,v]=q$, where $q=ih.$ Then the operators
$x=e^u$ и $y=e^v$ satisfy the relation $xy=qyx$.
The algebra, generated by $x$ and $y$, is called quantum torus.\\
In the preprint the author gives the definition of quantum algebraic torus
over the arbitrary field. Quantum algebraic tori can be characterized in
terms of exact sequences (Theorem 1).
The description of the center is given in Proposition 1. In the case of roots
of 1 QAT produces the fibering of central simple algebras on it's center
(Proposition 2).
QAT can be used for the construction of central simple algebras
and skew fields. It is unknown whether these
central simple algebras be crossed products or not.
The description of QAT of dimension 2 is given in terms of generating
elements and there relations. The case of roots of 1 is
considered (Proposition 3).\\
\\
\begin{center}
{\bf Quantum algebraic tori and exact sequences}
\end{center}
Let $L$ be a field and $Q=(q_{ij})_{ij=1}^n$ be a matrix with the entries
$q_{ij}\in L^*$, $q_{ij}q_{ji}=q_{ii}=1$. An algebra of twisted
Laurent polynomials is the algebra
$L_Q[x_1^{\pm 1},\ldots ,x_n^{\pm 1}]$ ,
generated by the elements $x_1^{\pm 1},\ldots ,x_n^{\pm 1}$
with the relations $x_ix_j=q_{ij}x_jx_i$.
The group of diagonal matrices
$D_n(L)=L^*\times \cdots \times L^*$ acts by automorphisms  $x_i\mapsto t_ix_i$
on $L_Q[x_1^{\pm 1},\ldots ,x_n^{\pm 1}]$. We shall call this
algebra an quantum form of diagonal torus $D_n$.\\
Let $T$ be an algebraic $k$-torus splited over the Galois extension $L/k$.
Denote by $\Gamma =Gal(L/k)$, $M={\hat T}$ the group of characters
 of the torus $T$,
definded over the algebraic closure ofthe field $k$. The group $\Gamma $
actes on the group $M$.
The map $T\mapsto M$ is an 1-1 correspondence
of the set of algebraic tori,
splited over $L$, to the set of finitely generated ${\bf Z}$-free
${\bf Z}[\Gamma ]$-modules. The both sets is defined up to equivalence[1,2].\\
Denote by $^aM$  the group $M$ in additive form.
We shall denote elements of the group $M$ by $t^m$, $m\in {^aM}$.
Choose the ${\bf Z}$-basis $e_1,\ldots ,e_n$ in $^aM$.
The charactors $t_i=t^{e_i}$ free generate $M$.
For $m=m_1e_1+\cdots +m_ne_n$ the element $t^m$ is equal to
$t_1^{m_1}\cdots t_n^{m_n}$.\\
Now we shall give the main definition.\\
{\bf Definition 1}. A quantum form of the $k$-torus $T$
is the $k$-algebra $A_Q$ with following properties:
1) $A_{Q,L}=A_Q\otimes _kL=L_Q[x_1^{\pm 1},\ldots ,x_n^{\pm 1}]$,
2) there exists the map of coaction
$R: A_Q\mapsto k[T]\otimes _kA_Q$ of the Hopf algebra $k[T]$ on $A_Q$,
such that
$R_L=R\otimes 1: A_{Q,L}\mapsto L[T]\otimes _LA_{Q,L} $
coinsides with the coaction of diagonal torus $R_L(x_i)=t_i\otimes x_i$.
We shall also call $A_Q$ quantum algebraic torus.\\
The group $T(k)$ actes by automorphisms on $A_Q$.\\
{\bf Definition 2}.
The homomorphism of $k$-algebras $\phi :A_Q\mapsto A'_{Q'}$
is called a homomorphism of quantum $k$-forms of algebraic torus $T$ is the
following diagram commutes:\\
$$\begin {array}{ccc}A_Q&\stackrel {R}{\longrightarrow }&k[T]\otimes A_Q\\
\downarrow \lefteqn {\phi }&& \downarrow \lefteqn {1\otimes \phi }\\
 A'_{Q'}&\stackrel {R'}{\longrightarrow }&k[T]\otimes A'_{Q'}\end {array}
\eqno (1)$$
Denote by $Quan_L(T)$ the category of quantum $k$-forms of
algebraic torus $T$, splited over$L$ and defineded up to an isomorphism.
Denote by $Ext_{\Gamma }^c(M,L^*)$ the category of central extensions
$\Gamma $-групп $$1\mapsto L^*\mapsto E\mapsto M\mapsto 1\eqno (2)$$
Here $L^*$ is in the center of $E$.
The elements of this category are also defineded up to an isomorphism.\\
{\bf Theorem 1}. The categories $Quan_L(T)$ and $Ext_{\Gamma }^c(M,L^*)$
are isomorphic.\\
{\bf Proof}.\\
1) Let $A_Q$ be a quantum form of the $k$-torus $T$. Denote by $E$
the group $A_{Q,L}^*$  of invertible elements of
the algebra $A_{Q,L}=L_Q[x_1^{\pm 1},\ldots ,x_n^{\pm 1}]$.
The group $E$ is generated by $x_1^{\pm 1},\ldots ,x_n^{\pm 1}, L^*$.
The group $E$ is a $\Gamma $-group with respect to the
natural action of the Galois group $\Gamma $ in $A_Q\otimes L=A_{Q,L}$.
We restrict $R_L$ on $E$.
From $R_L(x_i)=t_i\otimes
x_i$ we have $R_L:E\mapsto M\times E$. The map
$R_L$ s a homomorphism of $\Gamma $-groups.
Consider the projection $\pi _1:M\times E$ on $M$. The superposition
$\pi =\pi _1R_L:E\mapsto M $ is a homomorphism of  $\Gamma $-groups
with the kernel $L^*$.
There exists the exact sequence of $\Gamma $-groups
$1\mapsto L^*\mapsto E\stackrel {\pi }{\mapsto }M\mapsto1$.\\
The algebra $A_Q$ is uniquely reconstructed
from the $\Gamma $-group $E$
as the $k$-algebra of $\Gamma $-invariant elements in $L_Q[x_1^{\pm 1},\ldots ,
x_n^{\pm 1}]$.\\
2) The homomorphism of $k$-algebras $\phi :A_Q\mapsto A'_{Q'}$
is extended to the $\Gamma $-invariant $L$-homomorphism
$\phi \otimes 1:A_{Q,L}\mapsto A'_{Q',L}$.
We restrict it on $E$ and have the
$\Gamma $-homomorphismof groups
$\phi :E={(A_{Q,L})}^*\mapsto E'={(A'_{Q',L})}*$.
The commutativity of the following diagrams one can be derived from (1):
$$\begin {array}{ccc}E&\stackrel {R}{\longrightarrow} &M\times E\\
\downarrow \lefteqn {\phi }&& \downarrow \lefteqn {1\otimes \phi }\\
 E'&\stackrel {R'}{\longrightarrow} &M\times E'\end {array}\eqno (3)$$

$$\begin {array}{ccccccccc}1&\longrightarrow &L^*&\longrightarrow &
E&\stackrel {\pi}{\longrightarrow}&M&\longrightarrow &1\\
&&\downarrow \lefteqn {id}&& \downarrow \lefteqn {\phi}
&&\downarrow \lefteqn {id}&&\\
1&\longrightarrow &L^*&\longrightarrow &
E'&\stackrel {\pi '}{\longrightarrow}&M&\longrightarrow &1 \end{array}
\eqno (4)$$
Isomorphic quantum forms correspond to isomorphic exact sequences. \\
3) The aim of this section is to reconstructa quantum form from
the exact sequence. According to the section 1,
this reconstruction is unique.\\
Let $E$ be a central extension of the $\Gamma $-groups
$1\mapsto L^*\mapsto E\stackrel {\pi }{\mapsto} M\mapsto1$.
The group $L^*$ embedes in $E$.
Consider the $L$-algebra $L_Q[M]=E\otimes _{L^*}L$.
The algebra $L_Q[M]$ is an algebra of twisted Laurent polynomials.
The group $\Gamma $ acts in  $L_Q[M]$ by automorphisms.
The group $E$ is a group of invertible elements in $L_Q[M]$.\\
Consider the $k$-algebra $A_Q=L_Q[M]^{\Gamma }$.
We shall shall show that $A_Q$ is a quantum form of $T$.
This will finish the proof of the theorem.\\
3.1) The aim of this section is to proof $A_Q\otimes _kL=L_Q[M]$.
Choose the elememts $x_1,\ldots ,x_n$ из $E$ such, that
$\pi (x_i)=t_i$. Denote $x^m=x_1^{m_1}\cdots x_n^{m_n}$.
The elements
$x^m$, $m\in ^aM$ form the basis of $L_Q[M]$ over $L$.
The algebra $A_Q$ is $k$-spanned by the elements
$a=\alpha x^m+\alpha_1x^{m(1)}+\cdots +\alpha _sx^{m(s)}$ such, that
the set of degrees $m, m(1),\ldots ,m(s)$ is a $\Gamma $-orbit in $^aM$.
Let $H=St(m,\Gamma )$. Then for $h\in H$ we have $^hx^m=\gamma _hx^m$,
where $\gamma \in L^*$. The following equality holds $\gamma _{h_1h_2}=
\gamma _{h_1}{^{h_1}}\gamma _{h_2}$.
Therefore $\gamma _h$ is a cocycle
on $H$  with values in $L^*$. Apply the Hilbert Theorem 90 to the extension
$L$ over $F=L^H$.
There exists the element $\gamma \in L^*$ such,
that $\gamma _h={^h\gamma}\gamma ^{-1}$. Then $^h(\gamma ^{-1}x^m)=
\gamma ^{-1}x^m$.\\
Let $ e,g_1,\ldots , g_s$ be the representatives of cosets $G/H$ such, that
$g_i(m)=m(i)$.
The element $^{g_i}(x^m)$ is up to constant equales to $x^{m(i)}$.
Denote $x^{<m>}=\gamma ^{-1}x^m$, $x^{<m(i)>}={^{g_i}x^{<m>}}$.
The element $a$
equals to $\beta x^{<m>}+\beta _1x^{<m(1)>}+\cdots +\beta _sx^{<m(s)>}$,
where $\beta =\alpha \gamma \in F$ и $\beta _i={^{g_i}\beta }$.
Since the extension $L/k$ is separable, then the extension $F/k$
is also separable.
According to the Primitive Element Theorem there exists the element
$\delta \in F$ such, that $F=k(\delta )$.
If $f(x)$ is minimal $k$-polynomial of $\delta $, then $dim_kF=deg f(x)=s+1$.
Let $\delta _0=\delta,
\delta _1= {^{g_1}\delta }, \ldots , \delta _s= {^{g_s}\delta  }$ be all
roots of the polynomial $f(x)$.
Since $L/k$ is the Galois extension, then all roots are contained in $L$.
The algebra $A_Q$ is spanned over the field $k$ by the elements
$a_r=\delta _0^rx^{<m>}+\delta _1^rx^{<m(1)>}+\cdots +\delta _s^rx^{<m(s)>},
r=0,1,\ldots ,s$.
Denote by $W$ the matrix $\left (\delta _i^j\right )_{i,j=0}^s$
The determinant of $W$ is equals to $\prod _{i>j}(\delta _i-\delta _j)\ne 0$.
There exists the $s+1$-tuple $\nu =(\nu _0,\nu _1,\ldots ,\nu _s)$ such,that
$\nu W=(1,0,\ldots ,0)$. Then $\sum _{r=0}^s\nu_ra_r=x^{<m>}$  and
$x^m\in A_Q\otimes _kL$ for all $m\in {^aM}$. This proves
$A_Q\otimes _kL=L_Q[M]$.\\
3.2) The aim of this section is to construct the map of coaction on $A_Q$.
Recall, that the map $E\stackrel {\pi }{\mapsto}M$ is a homomorphism of
$\Gamma $-groups. Consider the $\Gamma $-homomorphism $R_L=(\pi ,id):
E\mapsto M\times E$. Extent it to the $\Gamma $-invariant
homomorphism $R_L:L_Q[M]\mapsto L_Q[T]\otimes _LL_Q[M]$.
The homomorphism $R_L$ coinsides with the coaction of diagonal torus
$R_L(x_i)=t_i\otimes x_i$.
Show, that the restriction $R_L$ on $A_Q$ be a coaction of $k[T]$.
The image  $R_L(A_Q)=R_L(L_Q[M]^{\Gamma })$ containes in
$(L[T]\otimes _LL_Q[M])^{\Gamma }$. Recall that $k[T]=L[T]^{\Gamma }$.
Let us show
$$(L[T]\otimes _LL_Q[M])^{\Gamma }=L[T]^{\Gamma }\otimes _kL_Q[M]^{\Gamma }$$
Let $c\in (L[T]\otimes _LL_Q[M])^{\Gamma }$. Consider the normal
basis $a_1,\ldots ,a_p$ in the extension $L/k$.
The  Galois group $\Gamma $  acts on the normal basis by permutatioms.
The element  $c$  can be represented in the form
$c=\sum _{i=1}^p(f_i\otimes \alpha _i)a_i$, where
$f_i,\alpha _i$ are $\Gamma $-invariant elements in $L[T]$, $L_Q[M]$.
For all $\sigma \in \Gamma$ we have $c={^\sigma c}$. Hence
$$\sum _{i=1}^p(f_i\otimes \alpha  _i)a_i=
\sum _{i=1}^p(f_i\otimes \alpha _i)({^{\sigma }a_i})=
\sum _{i=1}^p(f_i\otimes \alpha _i)a_{\sigma (i)}$$
Then $f_{\sigma (i)}\otimes \alpha _{\sigma (i)}=f_i\otimes a_i$
for all $i$.
The group $\Gamma $ transitively acts on the normal basis.
We have  $c=f_1\otimes \alpha _1(a_1+\cdots +a_p)\in L[T]^{\Gamma }
\otimes _kL_Q[M]^{\Gamma }$. $\Box $\\
\\
\\
\begin{center}
{\bf On construction of extensions}
\end{center}
Let $Q=(q_{ij})_{ij=1}^n$ be a matrix with the entries
$$q_{ij}\in L^*, q_{ij}q_{ji}=q_{ii}=1 \eqno (5)$$
Let $L_Q[M]=L_Q[x_1^{\pm 1},\ldots ,x_n^{\pm 1}]$ be an algebra
of twisted Laurent polynomials over the field $L$.
Denote by  $k_Q[T]$ or $A_Q$ the quantum form of torus $T$,
constructed by the extension (2).
This algebra coinsides with  $L_Q[M]^{\Gamma }$.
Consider bihomomorphism $Q:{^aM}\times {^aM}\mapsto L^*$, where
$Q(e_i,e_j)=q_{ij}.$
For $m=m_1e_1+\cdots +m_ne_n$ and $k=k_1e_1+\cdots +k_ne_n$ we have
$$Q(m,k)=\prod _{i,j=1}^n q_{ij}^{m_ik_j}$$
The bihomomorphism $Q$ satisfies the condition $Q(m,k)=Q(k,m)^{-1}.$
Denote $x^m=x_1^{m_1}\cdots x_n^{m_n}.$ We have
$$x^mx^k=Q(m,k)x^kx^m\eqno (6)$$
For all $\sigma \in \Gamma $ the  element ${^{\sigma}(x^m)}$ is equals
to $x^{\sigma m}$ up to constant:
$${^{\sigma}(x^m)}=\gamma _{\sigma}(m)x^{\sigma m},
\gamma _{\sigma }(m)\in L^*\eqno (7)$$
Acting by $\sigma $ on the left and right sides of the equation (6), we have
$$Q(\sigma m,\sigma k)={^{\sigma }Q(m,k)}\eqno (8)$$
As a result by the exact sequence (2) there constructed
the matrix $Q$ with the conditions (5) and (8).\\
{\bf Question }. Let $Q$ be a matrix with the conditions (5) and (8),
Can we reconstruct the exact sequence (2) ?
Or (the same in different terms) ,can we
define the action of Galois group on $L_Q[M]$ by automorphisms such that
(7) holds ?\\
The answer is positive in the following two cases.\\
1) Let $M$ is a permutation module. The group $\Gamma $ acts
by permutations on the basis
${^{\sigma }(e_i)}=e_{\sigma (i)}, i=\overline {1,n} $.
The equality (8) for the matrix $Q$ is equivalent to
$$q_{\sigma (i),\sigma (j)}={^{\sigma }q_{ij}}\eqno (9)$$
Let $L\{M\}=L\{x_1^{\pm 1},\ldots ,x_n^{\pm 1}\}$ be a free algebra
generated by $x_1^{\pm 1},\ldots ,x_n^{\pm 1}$.
The action ${^{\sigma }(x^m)}=x_{\sigma (i)}$ extents to the automorphism
of $L\{M\}$. Denote by $I$ the ideal, generated by the elements
$x_ix_j-q_{ij}x_jx_i$ in $L\{M\}$.
It follows from (9), that
the ideal $I$ is invariant under the action of the group $\Gamma $.
The group $\Gamma $ acts by autimorphisms in the
factor-algebra $L_Q[M]=L\{M\}/I.$\\
2) Let $M$ be an arbitrary module over the cyclic group
of the second order $\Gamma =\{e,\sigma : {\sigma }^2=1\}.$
The module $M$ is a direct sum indecomposable modules.
Indecomposable modules have rank 1 or 2.
In the first case the action $\sigma e_1=\pm e_1$.
In the second case $\sigma $ acts by permutations of basis elements
$\sigma e_1=e_2$, $\sigma e_2=e_1$.
the action $\Gamma $ extents to the action
in $L\{M\}$ by automorphisms:
${^{\sigma }x_1}=x_1^{\pm 1}$ for rank 1 and
${^{\sigma }x_1}=x_2$, ${^{\sigma }x_2}=x_1$ for rank 2.
Further the $\Gamma $-action in $L_Q[M]$ is consructed similar to 1).\\
\begin{center}
{\bf Center of $A_Q$}
\end{center}
Denote by $Z_Q[M]$ the center of $L_Q[M]$. The Galois group
$\Gamma $ acts in $L_Q[M]$ by automorphisms. The center $Z_Q[M]$
is invariant with respect to this action.
Denote ${^aM_c}=\{m\in {^aM}: Q(m,k)=1, \forall k\in {^aM}\}$.
It follows from (5) and (8) that the subgroup ${^aM_c}$
is invariant with respect to $\Gamma $-action.
Denote by $M_c$ the group ${^aM_c}$ in multiplicative form.
The center $Z_Q[M]$ is a linear space over the field
$L$, spanned by the elements $x^m, m\in {^aM_c}$.
We have $Z_Q[M]=L[M_c]$. \\
Denote by $Z_Q$ the center of $A_Q$.\\
{\bf Proposition 1}. The center $Z_Q$ coinsides with $Z_Q[M]^{\Gamma }$.\\
Proof. The subalgebra $Z_Q[M]^{\Gamma }$ consists of central elements
and is contained in $A_Q=L_Q[M]^{\Gamma }$.
Therefore $Z_Q[M]^{\Gamma }$
is contained in $Z_Q$. From the other side,
$Z_Q\subset Z_Q\otimes L\subset Z_Q[M]$.
Therefore $Z_Q=Z_Q^{\Gamma }\subset Z_Q[M]^{\Gamma }$.$\Box $\\
Denote by $T_c$ the $k$-torus corresponding to $\Gamma $-module $M_c$.
By the embedding of groups $i:M_c\mapsto M$
we can construct the embedding of $k$-algebras $i:k[T_c]\mapsto k[T]$
and the homomorphism $i_*:T\mapsto T_c$.\\
{\bf Corollary 1.1}. $Z_q=k[T_c]$. $\Box$ \\
Recall that the coaction $R: A_Q\mapsto k[T]\otimes _kA_Q$ lifts to
$R_L=R\otimes 1: A_{Q,L}\mapsto L[T]\otimes _LA_{Q,L} $ and coinsides
over $L$ with coaction of diagonal torus $R_L(x_i)=t_i\otimes x_i$. \\
Hence $R_L:Z_Q[M]\mapsto L[T_c]\otimes Z_Q[M]$.
The algebra $Z_Q[M]$ is spanned over the field $L$ by the elements
$x^m, m\in {^aM_c}$.
The algebra $L[T_c]$ is spanned over the field $L$ by the elements
$t^m, m\in {^aM_c}$.
The algebras $Z_Q[M]$ and $L[T_c]$ are isomorphic.
The coaction $R_L:L[T_c]\mapsto L[T_c]\otimes L[T_c]$ coinsides with
comultiplication on the diagonal torus $T_c\otimes L$.
The restriction $R_L$ on $A_Q$ is equal to $R$.
Similar to the proof of Theorem 1, one can show, that
$R:L[T_c]^{\Gamma }\mapsto L[T_c]^{\Gamma }\otimes _kL[T_c]^{\Gamma }$ .
The coaction $R$ coinsides with comultiplication
$k[T_c]\mapsto k[T_c]\otimes_kk[T_c]$\\
{\bf Corollary 1.2}. Coaction $R$ on the center $Z_Q=k[T_c]$ counsides with
comultiplication in $k[T_c]$.
The action of $T(k)$ on the center is a superposition
of the homomorphism $i_*$ and the regular action of
the torus $T_c$. $\Box$\\
Consider the case when $q_{ij}$ are all roots of degree $l$
of unity. Then $q_{ij}=\epsilon ^{s_{ij}}$, where  $s_{ij}\in {\bf Z}_l$,
$\epsilon \in L$ is a primitive root of degree $l$ of unity.
The center $Z_Q[M]$ of the algebra $L_Q[M]$ contains
the subalgebra $L[x_1^{\pm l},\ldots , x_n^{\pm l}]$.
The center $Z_Q$ of the algebra $A_Q$ contains the subalgebra
$Z_Q^{(l)}$= $L[x_1^{\pm l},\ldots , x_n^{\pm l}]^{\Gamma }$.
We shall call this algebra a $l$-center.
The algebra $Z_Q^{(l)}$ is an algebra of regular functions on the
$k$-torus $T^{(l)}$,which correponds to the
$\Gamma $-module ${^aM^{(l)}}$= $l({^aM})$= ${\bf Z}le_1+\cdots +{\bf Z}le_n$.
By the embedding of groups $i_l:M^{(l)}\mapsto M$ we can construct the
embedding of $k$-algebras $i_l:k[T^{(l)}]\mapsto k[T]$ and the homomorphism
$i_{l,*}:T\mapsto T^{(l)}$.
The homomorphisms $i_*$ and $i_{l,*}$ are isogenies.
Similar to the case of center the coaction on the $l$-center
coinsides with the comultiplication on $T^{(l)}$.
The action of $T$ on the $l$-center is a superposition
of the homomorphism $i_{l,*}$ and the regular action of the $T^{(l)}$. \\
\\
\\
\\
\begin{center}
{\bf Spectrum of quantum torus at roots of 1}
\end{center}
We shall consider the case
$q_{ij}=\epsilon ^{s_{ij}}$, where $s_{ij}\in {\bf Z}_l$,
$\epsilon \in L$ is a primitive root degree $l$ of unity.\\
Let $\chi \in T_c(k)=Hom_{k-alg.}(Z_Q,k)$. Denote by $\chi _l$
the restriction of $\chi $ on the $Z_Q^{(l)}$.
Extent $\chi $ to the $\Gamma $-invariant $L$-homomorphism
$\chi _L$ of the algebra $Z_Q[M]=Z_Q\otimes _kL$ into $L$.
Similar extent $\chi _l$ to the $\Gamma $-invariant $L$-homomorphism
$\chi _{l,L}$ of the algebra $Z_Q^{(l)}[M]=Z_Q^{(l)}\otimes _kL$ into $L$.\\
{\bf Notations 1}.\\
$P_{\chi ,L}=L_Q[M]/Ker {\chi _L}$,
$P_{\chi ,L}^{(l)}=L_Q[M]/Ker {\chi _{l,L}}$\\
$P_{\chi }=L_Q[M]/Ker {\chi }$,
$P_{\chi }^{(l)}=L_Q[M]/Ker {\chi _l}$\\
{\bf Notations  2}.\\
1) $ C_l(a,b,\omega )$  is an algebra over the field $L$
generated by the elements $x,y$ with the relations
$x^l=a$, $y_l=b$ $xy=\omega yx$, $a,b\in L^*$, $\omega = \epsilon ^k$,
$k$ divide $l$;\\
2) $C_l(a,b):=C_l(a,b,\omega )$ , if $\omega $ is a primitive root
of degree $l$ of 1; This algebra is called a cyclic algebra ;\\
3) $C_l(a):=L[x]/ (x^l-a)$, $a\in L^*$;\\
{\bf Definition}. The tensor product of some copies
of cyclic algebras is called a polycyclic algebra.\\
{\bf Prorosition 2}. \\
1)  $P_{\chi }^{(l)}$ and $P_{\chi,L}^{(l)}$ are semisimple algebras,
$P_{\chi}$ and $P_{\chi,L}$ are central simple algebras ,
$P_{\chi,L}$ is a polycyclic algebra;\\
2)$P_{\chi }^{(l)}\otimes _kL=P_{\chi ,L}^{(l)}$,
$P_{\chi}\otimes _kL=P_{\chi,L}$\\
{\bf Proof}. One can choose the set of generators in the
algebra $L_Q[M]$ of twisted Laurent polynomials as follows:
$x_i, y_i, z_j$, $i=\overline {1,p}$,
$j=\overline {1,n-2p}$ with the relations $x_iy_i=\epsilon ^{k_i}y_ix_i$,
$k_1\vert k_2,\cdots ,k_p\vert l $[3]. The else pairs of generators
commutes.
Denote $d_i=\frac {l}{k_i}$.\\
The center $Z_Q[M]$ of the algebra $L_Q[M]$ is generated by
$x_i^{d_i}$, $y_i^{d_i}$, $z_j$;
the $l$-center $Z_Q^{(l)}[M]$ is generated by $x_i^l$, $y_i^l$, $z_j^l$.
We have
$$P_{\chi,L}^{(l)}=\otimes _{i=1}^pC_l(a_i,b_i,\epsilon ^{k_i})\otimes _{j=1}
^{n-2p}C_l(c_j),$$
where $a_i=\chi (x_i^l)$, $b_i=\chi (y_i^l)$, $c_j=\chi (z_i^l)$.
The algebra $C_l(a,b,\omega )$ is a semisimple and is isomorphic
to the direct sum of cyclic algebras $C_d(\alpha , \beta )$,
where $\alpha $ and $\beta $ are contained in the extension of the field
$L$, $d$ divides $l$.
The algebras $P_{\chi ,L}^{(l)}$ and $P_{\chi }^{(l)}$ are semisimple.\\
The algebra $P_{\chi ,L}$ is isomorphic to the tensor product
$\otimes _{i=1}^pC_{d_i}(\alpha _i,\beta _i)$, where
 $\alpha _i=\chi (x^{d_i})$, $\beta _i=\chi (y^{d_i})$.
Therefore $P_{\chi ,L}$ is a polycyclic algebra and $P_{\chi ,L}$, $P_{\chi }$
are central simple algebras. This proves 1).\\
Note, that $(Ker \chi_L)^{\Gamma }=Ker \chi$ and
$Ker\chi \otimes L=Ker \chi _L$. Hence we have
$$P_{\chi }\otimes L=\left (A_Q/Ker \chi\right )\otimes _k L=
{A_Q\otimes L}/{Ker \chi \otimes L}={L_Q[M]}/{Ker \chi _L}=P_{\chi,L}$$
Similar one can get $P_{\chi }^{(l)}\otimes _kL=P_{\chi ,L}^{(l)}$. $\Box $\\
\\
\\
\begin{center}
{\bf Quantum algebraic tori of small dimension}
\end{center}
Let $L=k(\alpha )$, $\alpha ^2=D$, $D\in k$ be a quadratic
extension of the field $k$.
The Galois group of $L/k$ is $\Gamma =Gal(L/k)=\{1,\sigma :\sigma ^2=1\}$
In every $\Gamma $-module ${^aM}$ of rank 2 there exists
a basis such, that the matrix of $\sigma $ is equal to one of the following
matrices:
$\left (\begin{array}{cc} 1&0\\0&1\end{array}\right )$,
$\left (\begin{array}{cc} 1&0\\0&-1\end{array}\right )$,
$\left (\begin{array}{cc} -1&0\\0&-1\end{array}\right )$,
$\left (\begin{array}{cc} 0&1\\1&0\end{array}\right )$.
The algebra $L_Q[M]$ is generated by $x_1^{\pm 1},x_2^{\pm 1}$
with the relations $x_1x_2=qx_2x_1$.
The action of $\Gamma $ on ${^aM}$ extents to the action
on $L_Q[M]$ , if $q\in k$ in 1,3 cases and if ${^\sigma q}=q^{-1}$ (i.e.
$N(q)=1$) in 2,4 cases. Our aim to find the generators of
$A_Q=L_Q[M]^{\Gamma }$ and there relations.\\
1) $\sigma $ = $\left (\begin{array}{cc} 1&0\\0&1\end{array}\right )$.
The group $\Gamma $ trivially acts on ${^aM}$, $q\in k$.
The algebra $A_Q$ coinsides with the algebra of twisted
Laurent polynomials $k_Q[M]$.\\
In the case $q=1$ , we have $A_Q=A=k[x_1^{\pm 1},x_2^{\pm 1}]$ and
$T(k)=k^*\times k^*$\\
2) $\sigma $=$\left (\begin{array}{cc} 1&0\\0&-1\end{array}\right )$.
Denote $q=\lambda +\alpha \mu $. In this section
$N(q)=\lambda ^2-D\mu ^2=1$.
The algebra $A_Q$ is generated by
$x=x_1^{\pm 1}$, $z_1=\frac {x_2+x_2^{-1}}{2}$,
$z_2=\frac {x_2-x_2^{-1}}{2\alpha }$ with the relations
$$x\left (\begin{array}{c}z_1\\z_2\end{array}\right )=
\left (\begin{array}{cc}\lambda &D\mu \\
\mu &\lambda \end{array}\right )
\left (\begin{array}{c}z_1\\z_2\end{array} \right )x,
\quad z_1^2-Dz_2^2=1,\quad z_1z_2=z_2z_1$$
If $q=1$, then $A_Q=A$ coinsides with quotient algebra
$k[x_1^{\pm 1},z_1,z_2]$ with respect to the ideal generated by
$z_1^2-Dz_2^2=1$. The group
$T(k)=k^*\times T_{L/k}^{(1)}$, where $T_{L/k}^{(1)}=\{x\in L^* :N(x)=1\}$. \\
3) $\sigma $= $\left (\begin{array}{cc} -1&0\\0&-1\end{array}\right )$.
Denote
$y_1=\frac {x_1+x_1^{-1}}{2}$, $y_2=\frac {x_1-x_1^{-1}}{2\alpha }$,
$z_1=\frac {x_2+x_2^{-1}}{2}$, $z_2=\frac {x_2-x_2^{-1}}{2\alpha }$.
The algebra $A_Q$ is generated by $y_1,y_2,z_1,z_2$ with the relations
$y_1z_1\pm Dy_2z_2=q^{\pm }(z_1y_1\pm Dz_2y_2)$,
$y_2z_1\pm y_1z_2=q^{\pm }(z_2y_1\pm z_1y_2)$,
$ y_1^2-Dy_2^2= z_1^2-Dz_2^2= 1$, $[y_1,y_2]=[z_1,z_2]=0$.
In the case $q=1$ the algebra $A_Q=A$ coinsides with the quotient
algebra
$k[y_1,y_2,z_1,z_2]$ with respect to ideal generated by
$ y_1^2-Dy_2^2-1$, $ z_1^2-Dz_2^2-1$. The group $T(k)$ coinsides with
$T_{L/k}^{(1)}\times T_{L/k}^{(1)}$.\\
4)$\sigma $=$\left (\begin{array}{cc} 0&1\\1&0\end{array}\right )$,
$q=\lambda +\alpha \mu$.\\
Let $q\ne -1$.
Denote $u=\frac {x_1+x_2}{2}$, $v=\frac {x_1-x_2}{2\alpha }$,
$w=\frac {1+{^\sigma q}}{2}x_1x_2$. One can proof the following relations :
$$a) u^2-Dv^2=w,$$
$$b) uv-vu=-\frac {\mu}{1+\lambda}w,$$
$$c) w\left (\begin{array}{c}u\\v\end{array}\right )=
\left (\begin{array}{cc}\lambda &-D\mu \\
-\mu &\lambda \end{array}\right )
\left (\begin{array}{c}u\\v\end{array} \right )w$$
The formula $c)$ is derived from $a)$ and $b)$. \\
The algebra $A_Q$ is generated by $u,v, (u^2-Dv^2)^{-1}$
with the unique relation
$(1+\lambda )(uv-vu)+\mu(u^2-Dv^2)=0$.\\
In the case $q=-1$ denote
$u=\frac {x_1+x_2}{2}$, $v=\frac {x_1-x_2}{2\alpha }$,
$w=\alpha x_1x_2$.
One can proof the following relations :
$$a) u^2-Dv^2=0,$$
$$b) uv-vu=-\frac {w}{D},$$
$$c) w\left (\begin{array}{c}u\\v\end{array}\right )=
-\left (\begin{array}{c}u\\v\end{array} \right )w$$
The formula $c)$ is derived from $a)$ and $b)$. \\
The algebra $A_{-1}$ is generated by $u,v, (uv-vu)^{-1}$
with the unique relation
$u^2-Dv^2=0$.\\
In the case $q=1$ the algebra $A_Q=A=k[u,v,(u^2-Dv^2)^{-1}]$.
The group $T(k)$ coinsides with $L^*$.\\
{\bf Case of roots of 1}\\
Let $q$ be a primitive root of odd degree $l$ of unit.
In all of above cases 1)-4) the center $Z_Q[M]$ coinsides with the $l$-center.
$Z_Q^{(l)}[M]=L_Q[x_1^{\pm 1},x_2^{\pm 1}]$.
The center $Z_Q$ is a $k$-formof the center $Z_Q[M]$ (Proposition 1).
Our aim is to study the structure ofthe algebras $P(\chi )$ in common point.\\
Denote by $K=Fract(Z_Q)$ the field of fractions of $Z_Q$ and
$K_L=K\otimes _kL=Fract(Z_Q[M])$ is the field of fractions of $Z_Q[M]$.
Let $\chi :Z_Q\mapsto K$ be an embedding $Z_Q$ into $K$.
Extent $\chi $ to
$\Gamma $-invariant homomorphism $Z_Q[M]$ into $K_L$.
Denote $\chi (x_1)=a_1$, $\chi (x_2)=a_2$.
The algebra $P_{\chi ,L}$ counsides with $C_l(a_1,a_2)$ and  $P_{\chi }$
is it $K$-form.\\
{\bf Proposition 3}. In the cases 1),2),4) the skew field $P_{\chi }$
is a cyclic crossed product over the field $K$.
( The case 3) is unknown to the author).\\
{\bf Proof}. It is sufficient to prove that  $P_{\chi }$ contains the
maximal subfield, which is a Galois field with cyclic Galois group.
It is trivial in the case 1). In the cases 2) and 4)
there exists the element $y\in P_{\chi,L}$ such that ${^\sigma y}=y^{-1}$
( in the case 2) put $y=x_2$; in the case 3) put $y=x_1^{-1}x_2$ ).
In every case $y^l\in K$. As $q\in L$, then
$K_L(y)$ is a cyclic extension of the field $K_L$. Denote
$H=Gal(K_L(y)/K_L)= \{1,s,\ldots ,s^{l-1}: s^l=1\}$, ${^s}y=qy$.
Denote by $G$ the Galois group of the extension $K_L(y)/K$.
The group $G$ contains two subgroups $\Gamma $ and $H$.
Note that ${^{s\sigma }y}={^sy^{-1}}=q^{-1}y^{-1}=
{^{\sigma }(qy)}={^{\sigma s}y}$.
The subgroups $\Gamma $ and $H$ commutes. As $[K_L(y):K]=2l$, then
$G=H\times \Gamma$. The factor group $G/\Gamma $ is isomorphic to $H$.
The field $K_L(y)^{\Gamma }$ is a cyclic extension
of degree $l$ over $K$.
$\Box $\\
{\bf Bibliographi}: [1] A.Borel "Linear algebraic groups" New York-Amsterdam,
1969; [2] V.E.Voskresenskii "Algebraic tori", Moscow, Nauka, 1977;
[3] A.N.Panov "Skew fields of twisted rational functions and the skew field
of rational functions on $GL_q(n)$", St.Peterburg Math.J.Vol.7(1996),No1,
P.129-143.

\end{document}